\documentclass[10pt,twoside]{amsart}
\usepackage[english]{babel}
\usepackage[utf8]{inputenc}

\usepackage{multicol}

\usepackage{tabu}

\usepackage{amsmath}
\usepackage{amsthm}
\usepackage{amssymb}

\usepackage[colorlinks=true,linkcolor=blue,citecolor=blue,urlcolor=blue,citebordercolor={0 0 1},urlbordercolor={0 0 1},linkbordercolor={0 0 1}]{hyperref}
\usepackage{cite}

\numberwithin{equation}{section}

\newtheorem{thm}{Theorem}[section]

\newtheorem{obs}{Observation}
\newtheorem{ex}[thm]{Example}

\DeclareMathOperator{\Proj}{Proj}
\DeclareMathOperator{\D}{D}
\DeclareMathOperator{\Spec}{Spec}
\DeclareMathOperator{\colim}{colim}
\DeclareMathOperator{\QCoh}{QCoh}
\DeclareMathOperator{\Mod}{-Mod}
\DeclareMathOperator{\Hom}{Hom}
\DeclareMathOperator{\Cone}{Cone}
\DeclareMathOperator{\weight}{weight}

\newcommand{\bP}{\mathbb{P}}
\newcommand{\bA}{\mathbb{A}}
\newcommand{\cO}{\mathcal{O}}
\newcommand{\bG}{\mathbb{G}}
\newcommand{\bZ}{\mathbb{Z}}
\newcommand{\cC}{\mathcal{C}}
\newcommand{\cG}{\mathcal{G}}
\newcommand{\cU}{\mathcal{U}}
\newcommand{\bQ}{\mathbb{Q}}
\newcommand{\bR}{\mathbb{R}}
\newcommand{\bC}{\mathbb{C}}

\begin{document}

\markboth{\hfill{\rm Daniel Halpern-Leistner} \hfill}{\hfill {\rm On the structure of equivariant derived categories \hfill}}

\title{On the structure of equivariant derived categories}

\author{Daniel Halpern-Leistner}

\begin{abstract}
In this expository note, we discuss the results of \cite{Halpern-Leistner2015TheQuotient} on the structure of derived categories of equivariant coherent sheaves and the derived categories of geometric invariant theory quotients. We take a recent perspective, emphasizing the theory of restricted local cohomology. We also discuss several applications and concrete examples: studying the effects of birational modification on derived categories, constructing categorical completions of equivariant derived categories, and constructing actions of generalized braid groups on derived categories of GIT quotients. This is a contribution to the proceedings of the International Congress of Basic Science, held in July 2024.
\end{abstract}

\maketitle

\setcounter{tocdepth}{1}
\tableofcontents

\section{Introduction}

Throughout this note, we work over a base field $k$ of characteristic $0$. We will consider a reductive group $G$ acting on a finite type $k$-scheme $X$ that is projective over an affine scheme. We will assume that $X$ is equipped with a $G$-equivariant ample line bundle $L$. This is equivalent to saying that some positive power of $L$ comes from an equivariant closed embedding $X \hookrightarrow \bP^n \times \bA^m$, where $G$ acts linearly on the ambient space. 

In this setting, geometric invariant theory (GIT) provides an open $G$-equivariant subscheme $X^{\rm ss} \subset X$, defined as the union of the affine open subsets $\{\sigma \neq 0\}$ for all $\sigma \in \Gamma(X,L^n)^G$ for some $n \geq 0$. We will refer to the closed subset $S:= X \setminus X^{\rm ss}$ as the unstable locus. The main theorem of GIT is that $X^{\rm ss}$ admits a good quotient\cite{Seshadri1972QuotientGroups} by $G$
\[
X^{\rm ss} \to X^{\rm ss} /\!/ G = \Proj \left(\bigoplus_{n \geq 0} \Gamma(X,\cO_X(n))^G \right),
\]
known as the GIT quotient. In fact, $X^{\rm ss}$ only depends on the equivariant Neron-Severi class $\ell = c_1(L)$, and we will denote it $X^{\rm ss}(\ell)$ when we wish to make this dependence explicit.

Although it is common to regard GIT as a procedure taking $X \mapsto X/\!/G$, it is often more fruitful to regard it as two distinct steps: first one identifies the open substack of the quotient stack $X^{\rm ss} / G \subset X/G$, and second one constructs a good moduli space\footnote{This notion generalizes Seshadri's notion of good quotients. For a quotient stack, the two notions coincide.} $X^{\rm ss}/G \to X^{\rm ss}/\!/G$ for this open substack. In nice situations, $G$ will act with finite stabilizers on $X^{\rm ss}$, and the quotient map $X^{\rm ss}/G \to X^{\rm ss}/\!/G$ will be universally bijective. In general, the quotient map can identify many orbits in $X^{\rm ss}$. We will focus solely on the first step, i.e., we focus on the quotient stack $X^{\rm ss}/G$. 

Our goal is to survey some results comparing the derived category of equivariant coherent sheaves on $X$, or equivalently the derived category of coherent sheaves on the quotient stack $\D^{\rm b}(X/G)$, with that of the semistable locus $\D^{\rm b}(X^{\rm{ss}} / G)$. In order to do this, we will discuss a more general structure theorem for the unbounded derived category of equivariant quasi-coherent sheaves $\D^{\rm qc}(X/G)$. We will then discuss applications to comparing $\D^{\rm b}(X^{\rm ss}/G)$ for different GIT quotients of the same space $X$, to a categorification of the Atiyah-Bott localization formula, and to constructing exotic autoequivalences of $\D^{\rm b}(X^{\rm ss}/G)$.

This note is based on my Frontiers of Science Award lecture at the International Congress of Basic Science, held in July 2024 at the Beijing Institute for Mathematical Sciences and Applications. I thank the organizers for the award and for their hospitality.

\subsection{The origin and history}

The first general results comparing the geometry of the GIT quotient with the equivariant geometry of $X$ appear in the work of Kirwan \cite{Kirwan1985Cohomology31} and Atiyah-Bott in the infinite dimensional setting \cite{Atiyah1983}. These theorems state that the restriction map on equivariant cohomology $H^\ast_G(X;\bQ) \to H^\ast_G(X^{\rm ss};\bQ)$ admits a section, and hence is surjective. The theorems discussed in this note can be seen as a categorification of this result.

On the other hand, the idea that birational modifications of varieties should have predictable effects on derived categories of coherent sheaves goes back to the work of Bondal and Orlov \cite{Bondal1995SemiorthogonalVarieties}. A special case is the $D$-equivalence conjecture, which predicts that birational modifications that preserve the canonical bundle lead to equivalences of derived categories. Many of the ideas discussed below, and specifically the idea of using variation of GIT and ``window categories" to probe instances of this conjecture, originate in the work of Hori-Herbst-Page in physics \cite{Herbst2009B-typeVarieties}, and of Ed Segal in a mathematical context \cite{Segal2011EquivalencesB-Models}.

\section{Restricted local cohomology}

Local cohomology is, by definition, the object that measures the difference between the category of equivariant sheaves on $X$ and $X^{\rm ss}$. If $j : X^{\rm ss} \hookrightarrow X$ is the open immersion, then the local cohomology complex $R\Gamma_{S}(F)$ for any complex $F \in \D^{\rm qc}(X/G)$ can be defined by the exact triangle in $\D^{\rm qc}(X/G)$,
\begin{equation} \label{E:defining_triangle}
R\Gamma_S(F) \to F \to Rj_\ast(j^\ast(F)) \to.
\end{equation}

The second two terms in this triangle only depend on the underlying set of $S$, and not any particular choice of scheme structure, and therefore the same is true for $R\Gamma_{S}(F)$. The projection formula $Rj_\ast(j^\ast(F)) \cong F \otimes^L Rj_\ast(\cO_{X^{\rm ss}})$ shows that \eqref{E:defining_triangle} is obtained by tensoring the corresponding exact triangle for $\cO_X$ with $F$, and therefore $R\Gamma_{S}(F) \cong F \otimes^L R\Gamma_{S}(\cO_X)$.

\subsection*{The local cohomology of the unstable locus}
$S$ is cut out set-theoretically by a finite collection of $G$-invariant sections $\sigma_1,\ldots,\sigma_m \in \Gamma(X,L^n)^G$ for some $n$. Let $K(\sigma_1^j,\ldots,\sigma_m^j)$ denote the $G$-equivariant Koszul complex on $X$, whose $i^{\rm th}$ term is $K_i := \bigwedge^i (\cO_X(jn)^m)$ for $i=0,\ldots,m$. The differential $K_0 \to K_{1}$ is defined by $(\sigma_1^j,\ldots,\sigma_m^j) : \cO_X \to \cO_X(jn)^m$, and this is extended to a differential $K_i \to K_{i+1}$ for other $i$ using the graded Leibniz rule. There is a canonical homomorphism of complexes $K(\sigma_1^j,\ldots,\sigma_m^j) \to K(\sigma_1^{j+1},\ldots,\sigma_m^{j+1})$ induced on each term by the homomorphism $(\sigma_1,\ldots,\sigma_m) : \cO_X(jn)^m \to \cO_X((j+1)n)^m$. One has
\begin{equation} \label{E:explicit_computation}
R\Gamma_{S}(F) \cong \colim_j F \otimes K(\sigma_1^j,\ldots,\sigma_m^j).
\end{equation}

\subsection{Semiorthogonal decompositions}

The exact triangle \eqref{E:defining_triangle} is part of a certain categorical structure. If $\cC$ is a triangulated category, a \emph{semiorthogonal decomposition} $\cC = \langle \cC_1, \cC_2 \rangle$ consists of a pair of full triangulated subcategories such that $\Hom(E_2,E_1) = 0$ for and $E_1 \in \cC_1$ and $E_2 \in \cC_2$, and every $E \in \cC$ admits an exact triangle $E_2 \to E \to E_1 \to$. In this case the assignment $E \mapsto E_2$ is right adjoint to the inclusion functor $\cC_2 \subset \cC$, and $E \mapsto E_1$ is left adjoint to the inclusion functor $\cC_1 \subset \cC$. One useful fact is that given collections of objects $S_1, S_2 \subset \cC$ such that $R\Hom(E_2,E_1) \cong 0$ for any $E_1 \in S_1$ and $E_2 \in S_2$ and the subset $S_1 \cup S_2 \subset \cC$ generates $\cC$ as a triangulated category, one obtains a semiorthogonal decomposition $\cC = \langle \cC_1,\cC_2 \rangle$, where $\cC_i$ is the triangulated subcategory generated by $S_i$, for $i=1,2$.\footnote{This works with either meaning of generation: you can take it to mean the smallest triangulated subcategory containing a given subset, or to mean the smallest triangulated subcategory that is closed under direct summands.}

The local cohomology functor is the right projection functor for a semiorthogonal decomposition
\begin{equation}\label{E:local_cohomology_SOD}
    \D^{\rm qc}(X/G) = \langle \D^{\rm qc}(X^{\rm ss}/G) , \D^{qc}_S(X/G) \rangle,
\end{equation}
where $\D^{\rm qc}_S(X/G) \subset \D^{\rm qc}(X/G)$ denotes the full subcategory of complexes whose restriction to $X^{\rm ss}/G$ is acyclic, and we identify $\D^{\rm qc}(X^{\rm ss}/G)$ with a subcategory of $\D^{\rm qc}(X/G)$ via the derived pushforward $Rj_\ast$ along the open immersion, which is fully faithful.

The terminology of semiorthogonal decompositions is usually used for essentially small categories, but the local cohomology semiorthogonal decomposition \eqref{E:local_cohomology_SOD} only works for large categories of quasi-coherent sheaves. Although $\D^{\rm b}(X^{\rm ss}/G)$ is the Verdier quotient of $\D^{\rm b}(X/G)$ by the full subcategory $\D^{\rm b}_S(X/G)$ of complexes whose restriction to $X^{\rm ss}$ is acyclic, the local cohomology functor $R\Gamma_S(-)$ does not preserve bounded coherent complexes, and thus the inclusion $\D^{\rm b}_S(X/G) \subset \D^{\rm b}(X/G)$ does not admit a right adjoint.

In what follows, we will use the slightly more general notion of an $n$-term semiorthogonal decomposition $\cC = \langle \cC_1, \ldots, \cC_n \rangle$ of a triangulated differential graded (dg) category $\cC$. This is an ordered collection of full triangulated subcategories $\cC_i \subset \cC$ such that:
\begin{enumerate}
    \item $\Hom(E,F)=0$ for any $E \in \cC_j$ and $F \in \cC_i$ with $i<j$; and
    \item Every object in $E \in \cC$ admits a filtration $0=E_{n+1} \to E_n \to E_{n-1} \to \cdots \to E_1 = E$ such that $\Cone(E_{j+1} \to E_{j}) \in \cC_j$ for all $j=1,\ldots,n$.
\end{enumerate}
It can be shown that the filtration is unique up to natural isomorphism, and that the assignment $E \mapsto E_j$ is functorial in $E$. The functor $E \mapsto E_j$ is right adjoint to the inclusion of the full triangulated subcategory $\cC_{\geq j}$ generated by $\cC_j,\ldots,\cC_n$, and the functor $E \mapsto \Cone(E_j \to E)$ is left adjoint to the inclusion $\cC_{<j} \subset \cC$, where $\cC_{<j}$ is the triangulated subcategory generated by $\cC_1,\ldots,\cC_{j-1}$.

\subsection{The example of \texorpdfstring{$\bA^1/\bG_m$}{A1/Gm}}

Consider $X = \bA^1 = \Spec(k[x])$ with the action of $\bG_m$ given by $t \cdot z = t^{-1}z$. It is my opinion that nearly everything about geometric invariant theory can be understood by carefully studying this example!

A representation of $\bG_m$ is simply a graded vector space, and a $\bG_m$ action on $\bA^1$ is equivalent to a grading on its coordinate ring. In this case $x$ has degree $1$. $\cO_X(1)$ is the trivial line bundle, equipped with an action of $\bG_m$ such that its fiber (as a sheaf) at the point $0 \in \bA^1$ is a one dimensional vector space of degree $-1$.

The usual equivalence between the category of quasi-coherent sheaves on $\bA^1$ and $k[x]$-modules, taking $E \in \QCoh(\bA^1)$ to $\Gamma(\bA^1,E) \in k[x]\Mod$, extends to an equivalence between the category of equivariant quasi-coherent sheaves on $\bA^1$ and graded $k[x]$-modules. It will be more convenient for us to work with the category of graded modules and its derived category, rather than the categories of equivariant quasi-coherent sheaves.

The $\bG_m$-invariant global sections functor corresponds to taking the degree zero piece of the corresponding graded $k[x]$-module, and $\cO_X(1)$ corresponds to the free module $k[x] \cdot e$ with generator $e$ in degree $-1$. Therefore, $x^n$ can be regarded as an invariant section of $\cO_X(n)$ for any $n>0$, and the semistable locus with respect to $L=\cO_X(1)$ is $X^{\rm ss} = \bA^1 \setminus 0$. The unstable locus is $S = \{0\}$, cut out by the invariant section $x \in \Gamma(\cO_X(1))^{\bG_m}$.

The homomorphism $\cO_X \to j_\ast(\cO_{X^{\rm ss}})$ corresponds to the inclusion of graded modules $k[x] \hookrightarrow k[x^{\pm 1}]$, which shows that
\[
R\Gamma_S(M) \cong M \otimes^L R\Gamma_S(\cO_X) \cong M \otimes^L \left(k[x^{\pm 1}]/k[x][-1] \right).
\]
The graded $k[x]$ submodule of $k[x^{\pm 1}]/k[x]$ generated by $x^{-j}$ is finite dimensional, and we describe it as $k[x]/(x^j) \cdot e_{-j}$ with $e_{-j}$ a formal generator of degree $-j$. It has nonzero homogeneous elements of degrees $-1,\ldots,-j$. Writing the large module as a union $k[x^{\pm 1}]/k[x] = \bigcup_{j\geq 1} (k[x]/(x^j)) \cdot e_{-j}$, we have
\begin{equation}\label{E:specific_colimit}
R\Gamma_S(M) = \colim_j M \otimes^L \left( k[x]/(x^j) \cdot e_{-j}[-1] \right),
\end{equation}
which can be identified with \eqref{E:explicit_computation} by replacing the module $k[x]/(x^j) \cdot e_{-j}$ with its Koszul resolution, a complex in cohomological degree $0$ and $1$, $k[x] \to k[x] \cdot e_{-j}$, where $1 \mapsto x^j \cdot e_{-j}$. Note that for a free module $M$, the only non-vanishing cohomology group in \eqref{E:specific_colimit} is $H^1(R\Gamma_S(M))$.

\begin{obs}
    For a finitely generated free graded module $M$, $H^1(R\Gamma_S(M))$ is not finitely generated. However, for any $w \in \bZ$, the submodule spanned by elements of degree $\geq w$ for each term of the colimit \eqref{E:specific_colimit} is finite dimensional, and the colimit of these submodules stabilizes for $j \gg 0$. In particular, the submodule of $H^1(R\Gamma_S(M))$ spanned by elements of degree $\geq w$ is finite dimensional.
\end{obs}

Building on this observation, one can define a ``restricted" local cohomology functor $R\Gamma_{S}^{\geq w}(-)$ that assigns
\begin{equation}\label{E:restricted_local_cohomology_example}
R\Gamma_S^{\geq w} (\cO(n)) = \left\{ \begin{array}{ll} 0, &\text{if } n+w>-1 \\ k[x]/(x^{-n-w}) \cdot e_{n+w} (n)[-1] ,& \text{if } n+w \leq -1\end{array} \right. .
\end{equation}
The twist $(n)$ in the last expression denotes multiplying by a formal generator of degree $-n$, and $[-1]$ denotes a shift so that the module lies in cohomological degree $1$. To complete the definition of $R\Gamma_S^{\geq w}(-)$ as an exact functor on $\D^{\rm qc}(\bA^1/\bG_m)$, one could explicitly declare how the multiplication map $x : \cO(n) \to \cO(n+1)$ induces a homomorphism $R\Gamma_S^{\geq w} (\cO(n)) \to R\Gamma_S^{\geq w} (\cO(n+1))$ and verify that the functor respects isomorphisms of complexes of free modules. We will, however, take a different approach.

Each of the modules in \eqref{E:restricted_local_cohomology_example} can be filtered by powers of $x$, and the associated graded modules are of the form
\[
k(j) := k[x]/(x) \cdot e_{-j},
\]
where $e_{-j}$ is a formal generator of weight $-j$ with $w \leq -j \leq -n-1$. Let $\D^{\rm qc}_S(\bA^1/\bG_m)_{\geq w} \subset \D^{\rm qc}(\bA^1/\bG_m)$ denote the smallest full triangulated subcategory containing $k(j)$ for $j \leq -w$ and closed under arbitrary filtered colimits, and let $\D^{\rm b}_S(\bA^1/\bG_m)_{\geq w} \subset \D^{\rm qc}_S(\bA^1/\bG_m)_{\geq w}$ be the full subcategory of complexes with bounded coherent cohomology sheaves. Then we have the following:
\begin{thm}
    The inclusion $\D^{\rm qc}_S(\bA^1/\bG_m)_{\geq w} \subset \D^{\rm qc}(\bA^1/\bG_m)$ admits a right adjoint $R\Gamma_S^{\geq w}(-)$ that agrees with \eqref{E:restricted_local_cohomology_example}, and $R\Gamma_S^{\geq w}$ preserves $\D^{\rm b}(\bA^1/\bG_m)$, giving a right adjoint of the inclusion  $\D^{\rm b}_S(\bA^1/\bG_m)_{\geq w} \subset \D^{\rm b}(\bA^1/\bG_m)$.
\end{thm}
\begin{proof}[Sketch of proof]
    Because $\D^{\rm qc}(\bA^1/\bG_m)$ is the ind-completion of its category of compact objects $\D^{\rm b}(\bA^1/\bG_m)$, any exact functor $\D^{\rm b}(\bA^1/\bG_m) \to \D^{\rm b}_S(\bA^1/\bG_m)_{\geq w}$ extends uniquely to a functor $\D^{\rm qc}(\bA^1/\bG_m) \to \D^{\rm qc}_S(\bA^1/\bG_m)_{\geq w}$ that commutes with filtered colimits. If the first functor is right adjoint to the inclusion, then a formal argument shows that its extension is right adjoint to the inclusion of big categories. Therefore, all we must do is construct a right adjoint to the inclusion  $\D^{\rm b}_S(\bA^1/\bG_m)_{\geq w} \subset \D^{\rm b}(\bA^1/\bG_m)$.

    First, let $\D^{\rm b}_S(\bA^1/\bG_m)_{<w}$ denote the smallest full triangulated subcategory containing $k(j)$ for $j>-w$. One can compute that $R\Hom(k(j),k(i)) = 0$ for $j \leq -w < i$, and the objects $k(j)$ for $j \in \bZ$ generate $\D^{\rm b}_S(\bA^1/\bG_m)$ as a triangulated category. It follows that one has a semiorthogonal decomposition
    \[
    \D^{\rm b}_S(\bA^1/\bG_m) = \left\langle \D^{\rm b}_S(\bA^1/\bG_m)_{< w}, \D^{\rm b}_S(\bA^1/\bG_m)_{\geq w} \right\rangle.
    \]
    Let $\beta^{\geq w}$ and $\beta^{<w}$ denote the projection functors from $\D^{\rm b}_S(\bA^1/\bG_m)$ to the right and left semiorthogonal factor respectively.

    Next, we claim that for any $E \in \D^{\rm b}(\bA^1/\bG_m)$, the colimit \eqref{E:specific_colimit} stabilizes after applying $\beta^{\geq w}$, i.e., for all $j \gg 0$
    \[
    \begin{array}{rl} 0 &\cong \beta^{\geq w} \left(E \otimes^L \Cone \left(k[x]/(x^j) e_{-j} \xrightarrow{e_{-j} \mapsto x e_{-j-1}} k[x]/(x^{j+1})e_{-j-1} \right)[-1] \right) \\
    &\cong \beta^{\geq w}(E \otimes^L k(j+1))[-1]
    \end{array}
    \]
    This expression is functorial and exact in $E$, so it suffices to show this for $\cO(n)$ for any $n$, i.e., we must show that $\beta^{\geq w}(k(n+j+1))=0$ for all $j \gg 0$, which holds because for $j \gg 0$ the object $k(n+j+1) \in \D^{\rm b}_S(\bA^1/\bG_m)_{<w}$ by definition.

    We then define
    \[
    R\Gamma^{\geq w}_S(E) := \beta^{\geq w}\left(E \otimes^L k[x]/(x^j)e_{-j}[-1]\right),
    \]
    where $j \gg 0$ is chosen in the stable range for $E$. The fact that this functor is right adjoint to the inclusion $\D^{\rm b}_S(\bA^1/\bG_m)_{\geq w} \subset \D^{\rm b}_S(\bA^1/\bG_m)_{\geq w}$ is a formal consequence of the formula \eqref{E:specific_colimit} and the fact that $R\Hom(F, -)$ commutes with filtered colimits for $F \in \D^{\rm b}_S(\bA^1/\bG_m)_{\geq w}$.
\end{proof}

A completely analogous discussion allows one to construct a left adjoint $\beta^{<w}_S$ of the inclusion $\D^{\rm b}_S(\bA^1 / \bG_m)_{<w} \subset \D^{\rm b}(\bA^1/\bG_m)$. This time, the key fact one uses is that the pro-system $\{\beta^{<w}(E \otimes^L k[x]/(x^j))\}_{j \geq 0}$ stabilizes as $j \to \infty$. Putting all of this together gives the following structure theorems for $\D^{\rm b}(\bA^1/\bG_m)$ and $\D^{\rm qc}(\bA^1/\bG_m)$.

\begin{thm} \label{T:main_theta}
    There is a three-term semiorthogonal decomposition
    \[
    \D^{\rm b}(\bA^1/\bG_m) = \left\langle \D^{\rm b}_S(\bA^1/\bG_m)_{< w}, \cG_w,  \D^{\rm b}_S(\bA^1/\bG_m)_{\geq w} \right\rangle,
    \]
    where the restriction functor identifies $\cG_w \cong \D^{\rm b}((\bA^1)^{\rm ss}/\bG_m) \cong \D^{\rm b}(k\Mod)$. Replacing each of these categories with their closures under filtered colimits, one obtains an analogous semiorthogonal decomposition of $\D^{\rm qc}(\bA^1/\bG_m)$.
\end{thm}
\begin{proof}
    The existence of a such a decomposition is a formal consequence of the existence of the adjoint functors $R\Gamma_S^{\geq w}(-)$ and $\beta^{<w}_S(-)$. The category $\cG_w$ consists of those complexes $E$ for which $0 \cong R\Gamma_S^{\geq w}(E) \cong \beta^{<w}_S(E)$.
    
    However, in this case $\cG_w$ admits the more explicit description as the smallest full triangulated subcategory containing $\cO(-w)$. Taking this as the definition of $\cG_w$, the semiorthogonal decomposition follows from three observations:
    \begin{enumerate}
        \item $R\Hom(\cO(-w),k(j)) \cong k(j+w)^{\bG_m}$ vanishes when $j>-w$;
        \item $R\Hom(k(j), \cO(-w)) \cong k(-w-j+1)^{\bG_m}[-1]$ vanishes when $j \leq -w$; and
        \item using the short exact sequences $0 \to \cO(j-1) \to \cO(j) \to k(j) \to 0$ one can show that $\cO(-w)$ along with the objects $k(j)$ for $j \in \bZ$ generate $\D^b(\bA^1/\bG_m)$ as a triangulated category.
    \end{enumerate}  
    
    The fact that $\cG_w \cong \D^{\rm b}((\bA^1)^{\rm ss}/\bG_m)$ under restriction follows from the computation that $R\Hom(\cO(-w),\cO(-w)) \cong k$. It can also be deduced from a formal argument using the fact that $\D^{\rm b}_S(\bA^1/\bG_m)$ is generated by $\D^{\rm b}_S(\bA^1/\bG_m)_{\geq w}$ and $\D^{\rm b}_S(\bA^1/\bG_m)_{<w}$.
\end{proof}

\subsection{Restricted cohomology in general}

We now return to our initial setting, where $X$ is a $G$-scheme that is projective over an affine scheme. Theorem \ref{T:main_theta} is the template for a much more general phenomenon. We let $\D^-(X/G) \subset \D^{\rm qc}(X/G)$ denote the full subcategory of complexes with bounded above and coherent homology sheaves.


\begin{thm}\label{T:main_template}
    There is a ``restricted local cohomology" semiorthogonal decomposition $\D^{\rm qc}(X/G) = \left\langle \D^{\rm qc}(X/G)_{<\theta}, \D^{\rm qc}_S(X/G)_{\geq \theta} \right\rangle$. The restriction functor preserves $\D^-(X/G)$, leading to a three-term semiorthogonal decomposition
    \[
    \D^-(X/G) = \left\langle \cU_\theta^- , \cG_\theta, \cU_\theta^+ \right\rangle,
    \]
    where $\cU_\theta^+ = \D^-(X/G) \cap \D^{\rm qc}_S(X/G)_{\geq \theta}$, the outer two terms generate the subcategory $\D^{-}_S(X/G)$ of unstably supported complexes, and the restriction functor gives an equivalence $\cG_\theta \cong \D^{-}(X^{\rm ss}/G)$.
    
    Furthermore,if $X$ is smooth, then this decomposition restricts to a semiorthogonal decomposition of $\D^{\rm b}(X/G)$, and the middle category $\cG_\theta^{\rm b}:= \cG_\theta \cap \D^{\rm b}(X/G)$ is equivalent to $\D^{\rm b}(X^{\rm ss}/G)$ under restriction.
\end{thm}

Precisely specifying the categories in these decompositions lies outside the scope of this short survey. Instead, we will explicitly describe the categories as needed to understand the examples discussed below. The decomposition depends on a parameter $\theta$, valued in a real vector space, as we will see in the examples. 

A theorem of this kind is established in \cite[Thm.~2.10]{Halpern-Leistner2015TheQuotient} under more restrictive hypotheses. There, the categories $\cU^\pm_\theta$ and $\cG_\theta$ are described in terms of a canonical stratification of the unstable locus in geometric invariant theory. A much more general version of the theorem appears in \cite[Thm.~2.2.2]{Halpern-Leistner2021DerivedConjecture} -- it applies to any algebraic derived stack with a $\Theta$-stratification. The proof of these theorems follows similar steps to the proof of Theorem \ref{T:main_theta} sketched above.

\section{Application: wall crossing for derived categories}

For simplicity, we will describe this only in the case where $G = \bG_m$ acting on a smooth projective variety $X$. Let $X^{\bG_m} = Z_1 \bigsqcup \cdots \bigsqcup Z_N$ be the decomposition of the fixed locus into connected components. For each $i$, let $w_i = \weight(L|_{Z_i})$. Assume the indexing is such that $w_1 \leq \cdots \leq w_N$. For each $i$, we define $S_i^\pm = \{ x \in X | \lim_{t \to 0} t^{\pm 1} \cdot x \in Z_i\}$. 

Now fix a $G$-equivariant ample line bundle $L$, and for any $a \in \bZ$, let $L(a)$ be the different equivariant structure obtained by tensoring with the character of weight $-a$. The Hilbert-Mumford criterion for semistability with respect to $L(a)$ describes the unstable locus $S(a) \subset X$ as
\begin{equation}\label{E:unstable_locus}
    S(a) = \bigcup_{w_i<a} S_i^+ \cup \bigcup_{w_i>a} S_i^-,
\end{equation}
and the semistable locus is $X^{\rm ss}(a) = X \setminus S(a)$. The $S^\pm_i$ form a stratification of the unstable locus $S(a)$, in the sense that the closure of each $S^{\pm}_i$ lies in the union of $S^\pm_j$ for which $|w_j-a| \geq |w_i-a|$ and $w_j-a$ has the same sign as $w_i-a$.

As one varies $a$, the unstable locus \eqref{E:unstable_locus} changes only when $a=w_i$ for some $i$. Suppose $a$ has one of these critical values. Then for each $i$ with $a=w_i$, $S_i^+$ will be unstable for slightly larger values of $a$, $S_i^-$ will be unstable for slightly smaller values of $a$, and both $S_i^+$ and $S_i^-$ will lie in the semistable locus for $a$ itself. One can visualize this as unstable strata flipping from ascending to descending manifolds as $a$ varies. For $a$ sufficiently positive or negative, the unstable locus \eqref{E:unstable_locus} will be all of $X$, and the stratification of the unstable locus is known as the Bia{\l}ynicki-Birula stratification.

Let us describe the category $\cG^{\rm b}_w$ of Theorem \ref{T:main_template} more explicitly in this context. For each component $Z_i$ of $X^{\bG_m}$, let $T^*X|_{Z_i}^\pm$ be the direct summand of the cotangent bundle $T^*X|_{Z_i}$ whose weights are positive or negative respectively. We let $\eta_i^\pm = |\weight(\det(T^*X|_{Z_i}^\pm))|$, i.e, the sum of the positive or negative cotangent weights along $Z_i$. Likewise, $H^\ast(E|_{Z_i})$ splits as a direct sum of sheaves on which $\bG_m$-acts with weight $w$, and we say the weights of $H^\ast(E|_{Z_i})$ are the integers $w$ for which a non-trivial summand appears in some cohomology sheaf. Then we have
\[
\cG^{\rm b}_\theta(a) = \left\{ E \in \D^{\rm b}(X/G) \left|  \begin{array}{c} H^\ast(E|_{Z_i}) \text{ has weights in } \\
\left\{\begin{array}{ll} \theta_i + [0,\eta_i^+], &\text{if }w_i<a, \\ \theta_i + [0,\eta_i^-], & \text{if }w_i>a \end{array}\right. \end{array} \right. \right\},
\]
where the parameter $\theta = (\theta_1,\ldots,\theta_N) \in \bR^N$ has no integer coefficients. Theorem \ref{T:main_template} says that the restriction functor is an equivalence $\cG^{\rm b}_\theta(a) \cong \D^{\rm b}(X^{\rm ss}(a)/\bG_m)$ for all $a$.


As $a$ varies, the category $\cG^{\rm b}_\theta(a)$ only changes as $a$ crosses come critical value $a=w_i$. Let $\omega_X = \det(T^* X)$ denote the canonical bundle of $X$. If $\weight(\omega_X|_{Z_i})\geq 0$ for all $i$ with $a=w_i$, then $\eta_i^+ \geq \eta_i^-$, and so $\D^{\rm b}(X^{\rm ss}(a-\epsilon)/\bG_m) \cong \cG^{\rm b}_\theta(a-\epsilon) \subseteq \cG^{\rm b}_\theta(a+\epsilon) \cong \D^{\rm b}(X^{\rm ss}(a+\epsilon)/\bG_m)$. On the other hand, if $\weight(\omega_X|_{Z_i}) \leq 0$ for all $i$ with $a=w_i$, then $\D^{\rm b}(X^{\rm ss}(a+\epsilon)/\bG_m) \cong \cG^{\rm b}_\theta(a+\epsilon) \subseteq \cG^{\rm b}_\theta(a-\epsilon) \cong \D^{\rm b}(X^{\rm ss}(a-\epsilon)/\bG_m)$. Note that when $\weight(\omega_X|_{Z_i})=0$ for all $i$ with $a=w_i$, then one gets an equivalence of derived categories $\D^{\rm b}(X^{\rm ss}(a-\epsilon)/\bG_m) \cong \D^{\rm b}(X^{\rm ss}(a+\epsilon)/\bG_m)$

This paradigm, where there is a full category of $\D^{\rm b}(X/G)$ that is equivalent under the restriction functor to $\D^{\rm b}(X^{\rm ss}/G)$ for several different GIT quotients, has been fruitful for establishing new cases of the D-equivalence conjecture \cite{Halpern-Leistner2021DerivedConjecture}.

\section{Application: categorical localization}

One can also use Theorem \ref{T:main_template} to define the following ``completion" of the category $\D^{\rm b}(X/G)$ when $X$ is smooth and projective over an affine variety:
\[
\widehat{\D}(X/G) := \left\{F \in \D^{\rm qc}(X/G) \left| \begin{array}{c} H^i(F) \cong 0 \text{ for } |i|\gg 0, \\ F|_{X^{\rm ss}} \in \D^b(X^{\rm ss}/G), \text{ and}\\ \forall \theta, R\Gamma^{\geq \theta}_S(F) \in \D^{\rm b}(X/G) \end{array} \right. \right\}.
\]
Note that this category depends on the GIT parameter. In this section, we discuss some of the nice properties of this construction. All of these results are described in detail in \cite{Halpern-LeistnerAFormula}.
\begin{enumerate}
    \item $\widehat{\D}(X) \subset \D^{\rm qc}(X/G)$ is a triangulated, idempotent complete, symmetric monoidal subcategory (with respect to derived tensor product $\otimes^L$), so the Grothendieck group $K^0(\widehat{D}(X/G))$ has a canonical ring structure.\smallskip
    \item For any $F \in \widehat{D}(X/G)$, $H^\ast(X/G,F) \cong H^\ast(X,F)^G$ is finite dimensional, hence the index $\chi(X/G,F) := \sum_j (-1)^j \dim H^j(X,F)^G$ is well-defined.
\end{enumerate}
While $\widehat{\D}(X/G)$ contains certain infinite direct sums, its $K$-theory is not trivial, as is the case for $\D^{\rm qc}(X/G)$.

\begin{ex}
    If $X={\rm pt}$ and $G = \bG_m$, then $\D^{\rm qc}(X/G) = \D^{\rm qc}(B\bG_m)$ is the category of bounded complexes of graded vector spaces whose homology is finite dimensional in every weight and vanishes for all sufficiently large weights. In this case, $K^0(\widehat{\D}(\bG_m)) \cong \bZ(\!(u)\!)$, where $u$ is the class of the representation of weight $-1$.
\end{ex}

Let us return to the example of $X / \bG_m$ discussed in the previous section. To indicate the dependence of the completion on the GIT parameter $a$, we will write $\widehat{D}_a(X/\bG_m)$. We will abuse notation slightly to let $R\Gamma_{S_i^\pm} \cO_X$ denote the derived pushforward of the local cohomology complex from an equivariant open subscheme in which $S_i^\pm$ is closed. The stratification of the unstable locus \eqref{E:unstable_locus} by the $S_i^\pm$ induces a filtration of $\cO_X$ in $\D^{\rm qc}(X/\bG_m)$ whose associated graded objects are $Rj_\ast(\cO_{X^{\rm ss}})$ and $R\Gamma_{S_i^\pm} \cO_X \in \D^{\rm qc}(X/\bG_m)$ for the strata appearing in \eqref{E:unstable_locus}. One can use Theorem \ref{T:main_template} to show that all of these objects lie in $\widehat{\D}_a(X/\bG_m)$. We therefore get a decomposition of the unit in $K^0(\widehat{\D}_a(X/\bG_m))$,
\begin{equation} \label{E:decomp1}
    1 = [\cO_X] = [Rj_\ast(\cO_{X^{\rm ss}})] + \sum_{w_i<a} [R\Gamma_{S_i^+} \cO_X] + \sum_{w_i>a} [R\Gamma_{S_i^-} \cO_X]
\end{equation}
as a sum of mutually orthogonal idempotent elements.

For each $i$, we can define two completions of $\D^{\rm b}(Z_i/\bG_m) \cong \D^{\rm b}(Z_i \times (B\bG_m))$. The ``$+$" completion $\widehat{D}_+(Z_i/\bG_m)$ is the category of bounded complexes of graded quasi-coherent sheaves on $Z_i$ whose degree $w$ summand has coherent homology sheaves and is acyclic for $w \gg 0$. The ``$-$" completion $\widehat{D}_-(Z_i/\bG_m)$ has the same definition, except that the complex is acyclic in degree $w \ll 0$.

One can show that because the normal bundle $N_{Z_i} X$ of $Z_i$ splits as a sum of positive and negative weight bundles, with no summand in weight $0$, its Euler class
\[
e_\pm(N_{Z_i} X) := \sum_{i \geq 0} (-1)^i [ \bigwedge^i(N_{Z_i} X)^\ast],
\]
is a unit as an element of either completion $K^0(\widehat{\D}_+(Z_i/\bG_m))$ or $K^0(\widehat{\D}_-(Z_i/\bG_m))$. The subscript indicates which of these groups we regard $e_\pm(N_{Z_i} X)$ as an element in. Furthermore, for any $i$ with $w_i<a$, the pushforward functor along the closed immersion $\sigma_i : Z_i \hookrightarrow X$ maps $\widehat{\D}_+(Z_i/\bG_m)$ to $\widehat{D}(X(a)/\bG_m)$, and the same holds for $\widehat{\D}_-(Z_i/\bG_m)$ if $w_i>a$. Using various canonical filtrations of the local cohomology (see \cite{Halpern-LeistnerAFormula}) followed by the projection formula, one can show that for any $F \in \D^{\rm b}(X/\bG_m)$
\begin{equation}\label{E:reexpress_index}
\begin{array}{rl}
    \chi(X/\bG_m, R\Gamma_{S_i^\pm}(F)) &= \chi(Z_i/\bG_m, \sigma_i^\ast(F) / e_\pm(N_{Z_i} X))\\
    &= \chi(X/\bG_m, F \otimes^L (\sigma_i)_\ast( [\cO_{Z_i}] / e_\pm(N_{Z_i} X))).
\end{array}
\end{equation}

This leads to the following refinement of the decomposition \eqref{E:decomp1}: Let $K_{\rm nil} \subset K^0(\widehat{D}_a(X/G))$ be the subgroup of classes $\eta$ such that $\chi(X/\bG_m,V \otimes \eta) = 0$ for all $\bG_m$-equivariant vector bundles $V$ on $X$. We then define
\[
\widehat{K}^0_a(X/G) := K^0(\widehat{\D}_a(X/G)) / K_{\rm nil}.
\]
Interpreting \eqref{E:reexpress_index} as an equivalence modulo $K_{\rm nil}$, we arrive at the localization formula in $\widehat{K}^0_a(X/\bG_m)$
\begin{equation} \label{E:decomp2}
[\cO_X] = [Rj_\ast(\cO_{X^{\rm ss}(a)})] + \sum_{w_i<a} (\sigma_i)_\ast \left( \frac{[\cO_{Z_i}]}{e_+(N_{Z_i}X)} \right) + \sum_{w_i>a} (\sigma_i)_\ast \left( \frac{[\cO_{Z_i}]}{e_-(N_{Z_i}X)} \right).
\end{equation}

The formula \eqref{E:decomp2} is highly reminiscent of the Atiyah-Bott localization formula \cite{Atiyah1984TheCohomology}, except that it is in $K$-theory rather than cohomology, there is an additional term corresponding to the index over the semistable locus, and the terms coming from each fixed component $Z_i$ depend in a subtle way on the GIT parameter $a$.

This is a feature, not a bug. For $a \gg 0$ or $a \ll 0$, $X^{\rm ss}(a) = \emptyset$, and all of the fixed component terms have the same form. In this case, applying $\chi(X/G,F\otimes (-))$ to both sides of \eqref{E:decomp2} provides a formula for $\chi(X/\bG_m,F)$ as a sum of contributions from fixed components. On the other hand, for $a$ with $X^{\rm ss}(a) \neq \emptyset$, \eqref{E:decomp2} is a $K$-theoretic wall-crossing formula that can be used to express the difference $\chi(X^{\rm ss}(a_1)/\bG_m,F) - \chi(X^{\rm ss}(a_2)/\bG_m,F)$ as a sum over contributions from fixed components.

\section{Application: generalized braid group actions}

Let us return to the phenomenon discussed above, where applying Theorem \ref{T:main_template} for two or more different GIT quotients of $X$ results in the \emph{same} subcategory $\cG_\theta^{\rm b} \subset \D^{\rm b}(X/G)$. In this case the restriction functor
\begin{equation}\label{E:equivalence_principle}
    \cG_\theta^{\rm b} \xrightarrow{\cong} \D^{\rm b}(X^{\rm ss}(\ell_i)/G)
\end{equation}
is an equivalence for a sequence $\ell_1,\ldots,\ell_N$ of different GIT parameters for $X/G$. Here we focus on the parameter $\theta$, which is independent of $\ell$. Different values of $\theta$ give different derived equivalences $\D^{\rm b}(X^{\rm ss}(\ell_i)/G) \to \D^{\rm b}(X^{\rm ss}(\ell_j)/G)$, and hence autoequivalences of $\D^{\rm b}(X^{\rm ss}(\ell_i)/G)$.

In general, if $\ell$ is a GIT parameter such that the invertible sheaf $\omega_X \otimes_k \det(\mathfrak{g})$ is equivariantly trivializable on $X^{\rm ss}(\ell)$, where $\mathfrak{g}$ is the adjoint representation, then it is expected that $\D^{\rm b}(X^{\rm ss}(\ell')/G)$ are equivalent for all $\ell'$ that are close to $\ell$ (technically, in an adjacent GIT chamber) and such that $G$ acts with finite stabilizers on $X^{\rm ss}(\ell')$. We do not know how to verify \eqref{E:equivalence_principle} in this generality, but researchers have verified it in large classes of examples. For simplicity we will focus on certain examples where $X$ is a linear representation of $G$.

Let $T \subset G$ be a maximal torus, $W = N(T)/T$ the Weyl group, and $M = \Hom(T,\bG_m)$ the weight lattice of $T$. Any linear representation decomposes as a $T$-representation into a direct sum of characters, and we refer to the resulting set of elements of $M$ (with multiplicity) as the weights of the representation. Following \cite{Spenko2017Non-commutativeGroups}, we say that $X$ is quasi-symmetric if the sum of the weights of $X$ lying along any line in $M_\bQ$ vanishes. Under a certain genericity hypothesis on $X$, \cite{Halpern-Leistner2020CombinatorialEquivalences} defines a closed rational polytope $\nabla \subset M_\bQ$. For any $\theta \in M_\bQ^W$, we define $\mathcal{W}(\theta) \subset \D^{\rm b}(X/G)$ to be the full triangulated subcategory generated by the equivariant vector bundles $\cO_X \otimes U$ for representations $U$ whose weights lie in $\theta + \nabla$.

The mysterious coincidence is that even though one defines $\mathcal{W}(\theta)$ without reference to any GIT quotient, the categories coincide with the $\cG^{\rm b}_\theta(\ell)$ of Theorem \ref{T:main_template} for all sufficiently generic GIT parameters $\ell$. As a result, we have:

\begin{thm}\cite[Thm.~3.2]{Halpern-Leistner2020CombinatorialEquivalences}
    Let $X$ be a generic quasi-symmetric linear representation of $G$. There is a locally finite $M^W$-periodic hyperplane arrangement $\{H_\alpha\}_{\alpha \in I}$ in $M_\bQ^W$ such that the boundary $\partial (\theta + \nabla)$ contains no points of $M$ whenever $\theta \in M_\bQ^W \setminus \bigcup_\alpha H_\alpha$, and in this case the restriction functor $\mathcal{W}(\theta) \to \D^{\rm ss}(X^{\rm ss}(\ell)/G)$ is an equivalence for any $\ell$ that is not parallel to any of the $H_\alpha$.
\end{thm}

These equivalences of derived categories fit into the structure of a representation of the fundamental groupoid of a certain space, predicted by homological mirror symmetry. Specifically, we consider the complexification of the hyperplane arrangement $\{i H_\alpha \otimes \bC\}_{\alpha \in I}$ in $M_\bC^W$, and let $H = \bigcup_{\alpha \in I} iH_\alpha \otimes \bC \subset M_\bC^W$. Then we can assign the category $\mathcal{W}(\theta)$ to any point $\ell + i \theta \in M_\bC^W$ with $\ell=0$ and $\theta \notin H$, and we assign $\D^{\rm b}(X^{\rm ss}(\ell)/G)$ to any point $\ell + i \theta$ for which $\ell$ is not parallel to any $H_\alpha$. To the straight path $0 + i \theta \to \ell + i \theta$ with $\theta \notin H$ we assign the restriction functor $\mathcal{W}(\theta) \to \D^{\rm b}(X^{\rm ss}(\ell)/G)$, and to any path $\ell + i \theta_1 \to \ell + i \theta_2$ with $\ell$ not parallel to any $H_\alpha$ we assign the identity functor from $\D^{\rm b}(X^{\rm ss}(\ell)/G)$ to itself.

\begin{thm}\cite[Thm.~5.1]{Halpern-Leistner2020CombinatorialEquivalences}
    The above assignment of points of $M_\bC^W \setminus H$ to categories and paths in $M_\bC^W \setminus H$ to functors defines a functor from the fundamental groupoid $\Pi_1(M_\bC^W \setminus H)$ to the category of triangulated dg-categories and natural isomorphism classes of exact functors.
\end{thm}

In fact, the group of characters $M^W$ of $G$ acts on $\D^{\rm b}(X/G)$ by tensor product $F \mapsto \chi \otimes_k F$, and one can combine this action with the representation above to obtain a representation of the fundamental groupoid of the quotient $\Pi_1( (M_\bC^W \setminus H) / M^W )$.

\bibliography{references}{}
\bibliographystyle{plain}




\address{Department of Mathematics, University of ABC\\
 XXX City, USA.\\
\email{author1@abc.edu}}

\address{Department of Mathematics, University of ABC\\
 XXX City, China.\\
\email{author2@abc.edu.cn}}

\end{document}